\documentclass[11pt,a4,amssymb,epsfig]{article}

\usepackage{graphicx}
\usepackage{latexsym}
\usepackage{amsmath}
\usepackage{epsfig}
\usepackage{epic}
\usepackage{amssymb}

\newtheorem{lemma}{Lemma}[section]
\newtheorem{theorem}[lemma]{Theorem}
\newtheorem{proposition}[lemma]{Proposition}
\newtheorem{corollary}[lemma]{Corollary}

\newcommand{\pf}{\noindent{\em Proof: }}
\newcommand{\epf}{\hfill\hbox{\rule{3pt}{6pt}}\\}

\newcommand{\p}{{\mathcal P}}
\newcommand{\s}{{\Sigma}}
\newcommand{\T}{{\mathcal T}}
\newcommand{\Nc}{{\mathcal N}}
\newcommand{\A}{{\mathcal A}}

\newcommand{\cL}{{\mathcal L}}
\newcommand{\cP}{{\mathcal P}}

\renewcommand{\tilde}{\widetilde}

\begin{document}

\title{Phylogenetic Networks from Partial Trees}

\author{{\bf S.\,Gr{\"u}newald}\\
Department of Combinatorics and Geometry (DCG),\\
CAS-MPG Partner Institute for Computational Biology (PICB),\\
Shanghai Institutes for Biological Sciences (SIBS),\\
Chinese Academy of Sciences (CAS), China,\\
and \\
{\bf K.\,T.\,Huber} and {\bf Q.\,Wu}\\
School of Computing Sciences,\\
University of East Anglia,\\
Norwich, NR5 7TJ,\\
United Kingdom.\\
}
\date{\today}
\maketitle

\newpage

\noindent K.\,T.\,Huber\\
School of Computing Sciences,\\
University of East Anglia,\\
Norwich, NR5 7TJ,\\
United Kingdom.\\
email: katharina.huber@cmp.uea.ac.uk\\
FAX: +44 (0) 1603 593345\\

\noindent We would prefer to submit the final manuscript in LATEX.

\newpage

\begin{abstract}
A contemporary and fundamental problem faced by
many evolutionary biologists is how to puzzle together a collection $\cP$
of partial trees (leaf-labelled trees
whose leaves are bijectively labelled by species or, more
generally, taxa, each supported by e.\,g.\,a gene)
into an overall parental structure that displays
all trees in $\cP$. This already difficult problem is
complicated by the fact that the trees in $\cP$  regularly support
conflicting phylogenetic relationships and are not on the same
but only overlapping taxa sets. A desirable requirement on the
sought after parental structure therefore is that it can
 accommodate the observed conflicts. Phylogenetic networks are a popular tool
capable of doing precisely this. However, not much is
known about how to construct such networks from partial trees, a notable
exception being the $Z$-closure super-network approach and the
recently introduced $Q$-imputation approach.
Here, we propose the usage of closure rules to obtain such a network.
In particular, we introduce the novel $Y$-closure rule and show that
this rule on its own or in combination with one of Meacham's closure
rules (which we call the $M$-rule) has
some very desirable theoretical properties. In addition, we use the
$M$- and $Y$-rule to explore the
dependency of Rivera et al.'s ``ring of life''  on the fact that
the underpinning phylogenetic trees are all on the same data set.
Our analysis culminates in the presentation of 
a collection of induced subtrees from 
which this ring can be reconstructed. 
\end{abstract}

\section{Introduction}

Phylogenetic trees have proved an important tool for representing evolutionary
relationships. For a set $X$ of species (or, more generally, taxa)
these are formally defined as
leaf-labelled trees whose leaves are bijectively labelled by the elements of
$X$. Advances in DNA sequencing have resulted in ever more data on which
such trees may be based. Computational limitations however combined
with the need to understand species evolution
 have left biologists with the following fundamental problem
which we will refer
 to as {\em amalgamation problem}: given a collection
 $\cP$ of phylogenetic trees,
how can these trees  be amalgamated
into an overall parental structure that preserves the
phylogenetic relationships supported by
the trees in $\cP$?  The hope is that such a structure might help
shed light on the evolution of the underlying genomes
(and thus the species).
\begin{figure}[h]
\begin{center}
\input{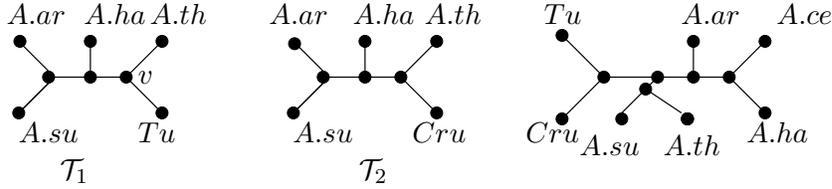}
\caption{\label{pete-trees}
3 phylogenetic trees which appeared in weighted
form in \cite{ML06} on subsets of the 7 plant species:
{\em A.thaliana (A.th)}, {\em A.suecia (A.su)}, {\em Turritis (Tu)}
{\em A.arenosa (A.ar)}, {\em A.cebennensis (A.ce)}, {\em Crucihimalaya (Cru)}
and {\em A.halleri (A.ha)}.
}
\end{center}
\end{figure}

In the ideal case that all trees in $\cP$
support the same phylogenetic relationships
(as is the case for trees  $\T_1$ and $\T_2$
depicted in Fig.~\ref{pete-trees})
this structure is known to be a phylogenetic tree and
a supertree method \cite{B04} may be used to reconstruct it.
For the above example the outcome $\T^*$ of such a method is $\T_1$
with species {\em Cru} (see Fig.~\ref{pete-trees} for full
species names) attached
via a pendant edge to the vertex labelled $v$.
It should be noted that $\T^*$
supports the same phylogenetic relationships
as $\T_1$ and $\T_2$ in the following sense:  For a finite set $X$, call
a bipartition $S=\{A,\tilde{A}\}$
of some subset $X'\subseteq X$ a {\em partial split} on $X$, or a
{\em partial (X)-split} for short, and
denote it by $A|\tilde{A}$ or, equivalently, by $\tilde{A}|A$
where $\tilde{A}:=X'-A$. In particular, call $S$ a
{\em (full) split} of $X$ if $X'=X$. Furthermore, say
that a partial $X$-split $S=A|\tilde{A}$ {\em extends} a partial
$X$-split $S'=B|\tilde{B}$ if either $B\subseteq A$ and
$\tilde{B}\subseteq \tilde{A}$
or $B\subseteq \tilde{A}$ and  $A\subseteq \tilde{B}$. Finally,
say that a phylogenetic tree
$\T$ {\em displays} a split $S=A|\tilde{A}$ if $S$ is a partial
split on the leaf set $\cL(\T)$ of $\T$ induced by deleting an edge of $\T$.
Then ``supports the same phylogenetic relationships''
means that for every split
 $S$ displayed by $\T_1$ or $T_2$ there
exists a split on $\cL(\T^*)$ that
extends $S$ and is displayed by $\T^*$.

\begin{figure}[h]
\begin{center}
\includegraphics[width=60mm,height=40mm]{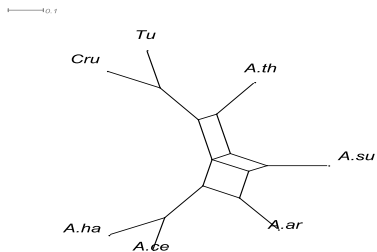}
\caption{\label{pete-network}
A circular phylogenetic network that represents all phylogenetic
relationships supported by the trees depicted in Fig.~\ref{pete-trees}
(see that figure for full species
names). }
\end{center}
\end{figure}

Due to complex evolutionary mechanisms such
as incomplete lineage sorting, recombination (in the
case of viruses), or lateral gene transfer (in case of bacteria) the
trees in $\cP$ may however support not
the same but conflicting phylogenetic relationships.
A phylogenetic network in the form of a split network (see
\cite{HM05,M05} for overviews) rather
than a phylogenetic tree is therefore the structure of choice
if one wishes to simultaneously represent all phylogenetic
relationships supported by the trees in $\cP$. An example in point
is the split network pictured in Fig.~\ref{pete-network}
which appeared as a weighted network in \cite{ML06}.
 With replacing ``edge'' in the definition of displaying
 by ``band of parallel edges'' and ``$\cL(\T)$'' by ``set of
network vertices of degree 1'' to obtain a definition for when a
 split network displays a split,
it is straight forward to check  that the network in Fig.~\ref{pete-network}
displays all splits displayed by the 3  trees pictured
in Fig.~\ref{pete-trees}.

It should be noted that phylogenetic networks such as the one depicted in
 Fig.~\ref{pete-network}
(see e.g. \cite{GEL03,HM06,HKLS05} for recently introduced other types of
phylogenetic networks) provide a means to visualize the complexity
of a data set and should not be thought of as
an explicit model of evolution. Awareness of this complexity does not
only  allow the exploration of a  data set but,
as is the case of e.g. hybridization networks \cite{HKLS05},
can also serve as starting
point for obtaining an explicit model of
evolution (see \cite{HB05} for more on this).

Apart from displaying all splits induced by the 3  trees depicted in
Fig.~\ref{pete-trees}, the network depicted in Fig.~\ref{pete-network}
has a further interesting feature. It is {\em circular}. In other words, if $X$
denotes the set of the 7 plant species under consideration,
then the elements of
$X$ can be arranged around a circle $C$ so that every split $S=A|\tilde{A}$ of
$X$ displayed by the
network can be obtained by intersecting $C$ with a straight line
so that the label set of one of the resulting 2
 connected components is $A$ and the label set of the other is $\tilde{A}$.

Although seemingly a very special type of phylogenetic network, circular
phylogenetic networks are a frequently used structure in phylogenetics
(see e.\,g.\,\cite{BM04,DH04,FCM06,GFDM06,H04})
as they do not only
naturally generalize the concept of a phylogenetic tree but are
also guaranteed
to be representable in the plane; a fact that greatly facilitates
drawing and thus analyzing them. However, although recently first steps
have been made
with regards to finding a solution to the
amalgamation problem in terms
of a phylogenetic network leading to the
attractive Z-closure \cite{HDKS04} and $Q$-imputation \cite{HGHM06}
approaches, very little is known about
a solution of this problem in terms of a circular phylogenetic network.

Intrigued by this and motivated by the fact that,
from a combinatorial point of view, phylogenetic trees
and networks are {\em split systems} (i.\,e.\,collections of full splits)
and that therefore the amalgamation problem boils down to
the problem of how to extend partial splits
on  some set  $X$ to splits on $X$,
we wondered whether closure rules for partial splits
could not be of help. Essentially mechanisms for splits' enlargement,
such rules have proved
useful for supertree construction and also underpin the above mentioned
$Z$-closure super-network approach. As it turns out, this is indeed the case.
As an immediate consequence of our main result (Corollary \ref{cyclosure}),
we obtain that for a collection  of partial splits that
can be ``displayed'' by a circular phylogenetic network $\Nc$,
the collection of (full) splits
generated by the closure rules in the centre of this paper
is guaranteed to be displayable by $\Nc$
and also independent
of the order in which the rules are applied. 

In a study aimed at shedding light into
the origin of eukaryotes, Rivera et al.
\cite{RL04} put forward the idea of a ``ring of life''
with the eukaryotic genome being the result of a fusion of
two diverse procaryotic genomes (see also \cite{ME04,RL04,SSSHRL05}). 
A natural and interesting question in this context is
how dependent Rivera et al.'s ring of life is on the 
fact that all underpinning trees are on the same taxa set. 
In the last section of this
paper, we provide a partial answer by presenting an example of
a collection of 
induced partial trees from which the ring of life can be reconstructed
using the $M$- and $Y$-rule.

The paper is organized as follows.
In Section \ref{closure-rules}, we first introduce some more terminology
and then restate one of
Meacham's closure rules (our $M$-rule) and introduce the novel $Y$-rule.
 In Section~\ref{relship}, we study the
relationship between the $M$- and $Y$-rule and
the closure rule that underpins the aforementioned $Z$-closure super-network
approach. In Section~\ref{properties}, we introduce the concept
of a circular collection of partial splits and show that both the
$Y$- and  $M$-rule preserve circularity (Proposition~\ref{circ}).
In Section~\ref{spl-sequence}, we introduce the concept of
a split closure and show that for certain
collections of partial splits this closure is independent
of the order in which the $Y$-rule and/or $M$-rule are/is applied
(Theorem~\ref{split-closure}). This result lies at the heart of
Corollary~\ref{cyclosure}. In Section \ref{rivera-exmple}, we explore 
the dependency of Rivera et al.'s ring of life  on the fact that
the underpinning trees are all on the same data set

Throughout the paper, $X$ denotes a finite set and the terminology and
notation largely follows \cite{SS01}.

\section{Closure rules
}\label{closure-rules}

We start this section by introducing some
additional terminology and notation. Subsequent to this, we first
restate Meacham's rule (which we call the $M$-rule) and then
introduce a novel closure rule which we call the $Y$-rule.

Let $\s(X)$ denote the collection of all partial splits of $X$
and suppose $\s\subseteq \s(X)$. Then a partial
split $S\in \s$
that can be extended by a partial split
$S'\in \s-S$ is called {\em redundant}.
The set obtained by
removing redundant elements from $\s$ is denoted by $\s^-$.
If $\s=\s^-$ then $\s$ is called {\em irreducible} and the set of
all irreducible subsets in $\s(X)$ is denoted by $\p(X)$.
Note that the relation ``$\preceq$'' defined for any two
(partial) split collections $\s,\s'\in \p(X)$ by putting
$\s\preceq \s'$ if every partial split in
$\s$ is extended by a partial split in $\s'$ is a partial order on
$\p(X)$ \cite{SS01}.

Suppose for the following that $\theta$ is a closure rule, that is, 
a replacement rule that replaces a collection $\A\subseteq \Sigma(X)$ of
partial splits that satisfy some condition $C_{\theta}$
by a collection $\theta(\A)\subseteq \Sigma(X)$
whose elements are generated in some
systematic way from the partial splits in $\A$ (see e.g. the $M$-
and the $Y$-closure rules presented below for two such systematic ways).
Suppose $\s,\s'\in \p(X)$ are two irreducible collections of
partial splits and $C_{\theta}(\Sigma)$ is the set of all subsets
of $\Sigma$ that satisfy $C_{\theta}$. If there exists
some subset $\A\in C_{\theta}(\Sigma)$ such that
$\s'=(\s\cup \theta(\A))^-$ then we say that $\s'$ is {\em obtained
from $\s$ via a single application of $\theta$}.
Finally, if for every subset $\A\in C_{\theta}(\Sigma)$ we have 
$\theta(\A)^-\preceq \s$ then
we call an application of $\theta$ to $\s$ {\em trivial} and 
say that $\s$ is {\em closed}
with respect to $\theta$.


We are now in the position to present
the 2 closure rules we are mostly concerned with in this paper:
the $M$-rule which is originally due to Meacham \cite{M83}
and the novel $Y$-rule. We start with Meacham's rule.

\subsection{The $M$-rule}\label{rule-(TR)}
Suppose $S_1,S_2\in \s(X)$ are two distinct partial splits of $X$.
Then the {\em $M$-rule} $\theta_M$ is as follows:
\begin{enumerate}
\item[($\theta_M$)] If there exists some $A_i\in S_i$, $i=1,2$ such that
\begin{eqnarray}
A_1\cap A_2\not=\emptyset\mbox{ and }
\tilde{A_1}\cap \tilde{A_2}\not=\emptyset
\end{eqnarray}
then replace $\A=\{S_1,S_2\}$ by the set $\theta^{\{A_1,A_2\}}_{M}(\A)$
which comprises of $\A$ and,
in addition, also the partial splits
$$
S_1'=(A_1\cap A_2)|(\tilde{A_1}\cup \tilde{A_2}) \mbox{ and }
S_2'=(\tilde{A_1}\cap \tilde{A_2})|(A_1\cup A_2).
$$
\end{enumerate}
In case the partial splits $S_1$ and $S_2$ are such
that there is no ambiguity with regards to the identity of the sets $A_1$
and $A_2$ in the statement of the $M$-rule or they are irrelevant to
the discussion, we will simplify
 $\theta^{\{A_1,A_2\}}_M(\A)$ to $\theta_{M}(\A)$. Clearly, such ambiguity
cannot arise if $S_1$ and $S_2$ are {\em compatible}, that is, 
there exist subsets $D_i\in S_i$, $i=1,2$ such that $D_1\cap D_2=\emptyset$.
However if $S_1$ and $S_2$ are {\em incompatible}, that is, not compatible
then caution is required.

Note that if $A_1$ and $A_2$ as in the statement of the $M$-rule
are such that $A_2\subseteq A_1$ and
$\tilde{A_1}\subseteq \tilde{A_2}$, then it is easy to verify that
$\theta_{M}$
 applies trivially to $\A$. 
Also note that for any $\s\in \p(X)$ and any two distinct
partial splits $S_1,S_2\in \s$, we have
$$
\s\preceq (\s\cup\theta_{M}(\{S_1,S_2\}))^-.
$$
Finally, note that any phylogenetic tree  on $X$
that displays the partial splits in some set $\s\in \cP(X)$  also displays
the partial splits in $(\s\cup\theta_{M}(\{S_1,S_2\}))^-$, $S_1,S_2\in \s$.


\subsection{The $Y$-rule}\label{rule-(SK)}
\begin{figure}[h]
\begin{center}
\input{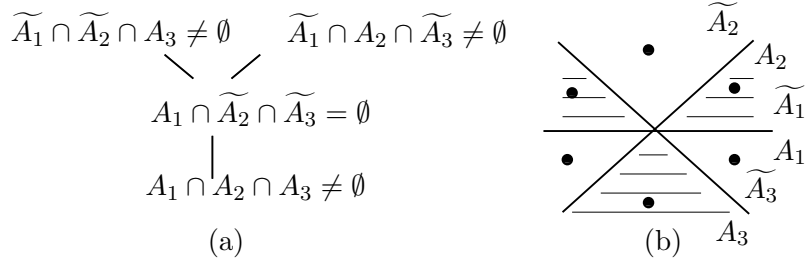}
\caption{\label{C_SK}
(a) A graphical representation of Condition (\ref{SK-rule}) in
the form of a $Y$.
(b) An example of three splits, depicted in bold lines, that
satisfy Condition (\ref{SK-rule})  -- see text for details.
}
\end{center}
\end{figure}

Suppose $S_i\in \s(X)$,
$i=1,2,3$, are three distinct partial splits of $X$.
Then the {\em $Y$-rule} $\theta_Y$ is as follows:
\begin{enumerate}
\item[($\theta_Y$)] If there exists some $A_i\in S_i$, $i=1,2,3$ such that
\begin{eqnarray}\label{SK-rule}
&&\emptyset\not \in\{ A_1\cap A_2\cap A_3,
\tilde{A_1}\cap \tilde{A_2}\cap A_3,
 \tilde{A_1}\cap A_2\cap \tilde{A_3}\}
\mbox{ and } \nonumber \\
&&A_1\cap \tilde{A_2}\cap\tilde{A_3} =\emptyset.
\end{eqnarray}
(see Fig.~\ref{C_SK}(a) for a graphical interpretation), 
then replace $\A=\{S_1,S_2,S_3\}$
by the set $\theta_Y^{\{A_1,A_2,A_3\}}(\A)$ which
comprises of the partial splits
\begin{eqnarray*}
S_1'&=&\tilde{A_1}\cup (\tilde{A_2}\cap \tilde{A_3})|A_1,
S_2'=A_2\cup (A_1\cap \tilde{A_3})|\tilde{A_2}, \mbox{ and }\\
S_3'&=&A_3\cup (A_1\cap \tilde{A_2})|\tilde{A_3}.
\end{eqnarray*}
\end{enumerate}

Although the condition in (\ref{SK-rule}) might look quite strange at first
sight, the class of triplets of partial splits that satisfy it is very
rich. For example, suppose that $S_i=A_i|\tilde{A_i}$, $i=1,2,3$ are splits of
$X$ that can be arranged in the plane as indicated in Fig.~\ref{C_SK}(b)
where
each bold, straight line represents one of $S_i$, $i=1,2,3$
and the dots represent non-empty triplewise intersections of
the parts of $S_i$, $i=1,2,3$,
in which they lie. For example,
the dot in the bottom wedge represents the
intersection $A_1\cap A_2\cap A_3$.
The shaded regions correspond
to the 3 non-empty intersections mentioned  in the statement of the
$Y$-rule. The partial
splits $S_i'=A_i'|\tilde{A_i'}$ $i=1,2,3$ obtained by restricting
$S_1$, $S_2$, and $S_3$ to different subsets of $X$ so that the
shaded regions remain non-empty form a triplet of partial splits that
satisfy (\ref{SK-rule}).

As the example of set $\A$ comprising the three partial splits
$S_1=145|2367$, $S_2=1357|246$, and $S_3=127|356$ shows different
choices of the sets $A_i$, $i=1,2,3$ lead to different sets
$\theta_Y^{\{A_1,A_2,A_3\}}(\A)$. For example, if
$A_1:=\{1,4,5\}$, $A_2:=\{1,3,5,7\}$, and $A_3:=\{1,2,7\}$
then (\ref{SK-rule}) is satisfied and
$\theta_Y^{\{A_1,A_2,A_3\}}(\A)= \{S_1,S_2,1247|356\}$.
If however $A_1$ and $A_2$ are as before and $A_3:=\{3,5,6\}$, 
then (\ref{SK-rule}) is also satisfied and
$\theta_Y^{\{A_1,A_2,A_3\}}(\A)$ is the set $ \{S_1,S_2,127|3456\}$.
Following our practise for the $M$-rule, for $\A=\{S_1,S_2,S_3\}$ we simplify
$\theta^{\{A_1,A_2,A_3\}}_{Y}(\A)$ to $\theta_{Y}(\A)$
if the  partial splits
$S_i$, $i=1,2,3$ are such
that there is no ambiguity with regards to the identity of
the sets $A_i$, $i=1,2,3$,
in the statement of the $Y$-rule or they are irrelevant to the discussion.

Note that
if $A_i$, $i=1,2,3$ as in the statement of the $Y$-rule are such that,
in addition,
$\emptyset\not=A_1\cap \tilde{A_2}\subseteq  A_3$, 
$\emptyset\not=A_1\cap \tilde{A_3}\subseteq  A_2$,
 $\emptyset\not=\tilde{A_3}\cap \tilde{A_2}\subseteq  \tilde{A_1}$,
and $\tilde{A_1}\cap A_2\cap A_3\not=\emptyset$  it is easy to see that 
$\theta_Y$ applies trivially to $\A$. 
Also note that
for any $\s\in \p(X)$ and any 3 partial splits $S_1,S_2,S_3\in\Sigma$
of $X$, we have
$$
\s\preceq (\s\cup\theta_Y(\{S_1,S_2,S_3\})^-.
$$



\section{First closure rule relationships}\label{relship}

In this section we first restate the {\em $Z$-(closure) rule}  which
was used in \cite{HDKS04} in the context
of a supernetwork construction approach
and then investigate the
relationship between the $Y$-, $M$-, and $Z$-rule.

Also originally
due to Meacham \cite{M83}, the $Z$-rule $\theta_Z$ can be restated as follows:
Suppose $S_1,S_2\in \s(X)$ are two distinct partial splits of $X$.
\begin{enumerate}
\item[($\theta_Z$)] If there exists some $A_i\in S_i$, $i=1,2$ such that
\begin{eqnarray}\label{condition1}
\emptyset \not\in\{ A_1\cap A_2, A_2\cap \tilde{A_1},
\tilde{A_1}\cap \tilde{A_2}\} \mbox{ and } A_1\cap \tilde{A_2}=\emptyset
\end{eqnarray}
then replace $\A=\{S_1,S_2\}$ by the set $\theta_{Y}(\A)$
which comprises of the partial splits
$(\tilde{A_1}\cup \tilde{A_2})|A_1$ and $ \tilde{A_2}|(A_1\cup A_2)$.
\end{enumerate}
Note that any two compatible partial splits of $X$ satisfy the
condition in (\ref{condition1}).

With this third closure rule at hand we are now in the
position to present a first easy to verify result. Suppose
$S_1$, $S_2$, and $S_3$ are 3 distinct partial splits of $X$
such that there exist parts
$A_i\in S_i$, $i=1,2,3$ as in the statement of the $Y$-rule. If, in addition,
$\tilde{A_1}\cap\tilde{A_2}\cap \tilde{A_3}\not=\emptyset$ and
$A_1\subseteq A_2\cup A_3$,  then
the partial split $A_1|\tilde{A_1}\cup(\tilde{A_2}\cap\tilde{A_3})$
generated by $\theta_Y$
is also generated by first applying $\theta_M$ to
$S_2$ and $S_3$ (with regards to $\tilde{A_2}\cap \tilde{A_3}\not=\emptyset$)
and then applying $\theta_M$ to the resulting partial split
$A_2\cup A_3|\tilde{A_2}\cap\tilde{A_3} $ and $S_1$. 

In addition, we have the following result whose
straight forward proof we leave to the reader.

\begin{proposition}
Suppose $S_1$ is a full split of $X$. Then the following
statements hold.
\begin{enumerate}
\item[(i)] If $S_2$ is a partial $X$-split and
$\theta_Z$ applies to $\Sigma=\{S_1,S_2\}$,
then
$$
\theta_{M}(\Sigma)^-=\theta_{Z}(\Sigma).
$$
\item[(ii)] If $S_2$ and  $S_3$
are partial $X$-splits so that $\theta_Y$
applies to $\Sigma=\{S_1,S_2,S_3\}$ and $\theta_Z$ applies to
$\{S_1,S_2\}$ and $\{S_1,S_3\}$.
Then
$$
(\theta_{Y}(\Sigma)\cup \bigcup_{j=2,3}\theta_{M}(S_1,S_j))^-
=(\bigcup_{i\in\{2,3\}}\theta_{Z}(S_1,S_i))^-.
$$
\end{enumerate}
\end{proposition}

\section{Closure rules and  weakly compatible
collections of partial splits}\label{properties}

In this section we introduce the notion of a
 weakly compatible collection of partial splits and
study properties of the $Y$- and $M$-rules
regarding such collections. A particular focus lies
on the study of circular collections of partial splits which we
also introduce. As we will see,
they form a very rich subclass of such collections of partial splits.

\subsection{Weakly compatible collections of partial splits}

We start this section with a definition that generalizes the
concept of weak compatibility for (full) splits of $X$ \cite{BD92} to partial
splits of $X$. Suppose  $S_i=A_i|\tilde{A_i}\in\s(X)$, $i=1,2,3$,
are three partial $X$-splits. Then we call $S_1,S_2,S_3$
{\em weakly
compatible} if  at least one of the
four intersections
\begin{eqnarray}\label{thinly}
 A_1\cap A_2\cap A_3, \tilde{A_1}\cap \tilde{A_2}\cap A_3,
 \tilde{A_1}\cap A_2\cap \tilde{A_3}, A_1\cap \tilde{A_2}\cap\tilde{A_3}
\end{eqnarray}
is empty\footnote{In the definition of
weak compatibility for full splits, $S_1$, $S_2$ and $S_3$ are full
splits and the condition in (\ref{thinly}) is the same (see \cite{BD92}).}.
Since the roles of $A_i$ and $\tilde{A_i}$ in $S_i$, $i=1,2,3$,
 can be interchanged without changing $S_1,S_2,S_3$
we have that $S_1,S_2,S_3$ are weakly compatible if and
only if at least one of the four intersections
\begin{eqnarray*}
 \tilde{A_1}\cap \tilde{A_2}\cap \tilde{A_3}, A_1\cap A_2\cap \tilde{A_3},
 A_1\cap \tilde{A_2}\cap A_3, \tilde{A_1}\cap A_2\cap A_3
\end{eqnarray*}
is empty.
More generally, we call a collection
$\s \subseteq \s(X)$ of partial $X$-splits {\em  weakly
compatible} if every three partial splits in
$\s$ are  weakly compatible.  To give an
 example, the partial splits
$S_1=123|4567$, $S_2=124|3567$, and $S_3=235|146$
are  weakly compatible whereas the partial splits
$S_3$, $S_4=24|135$, and $S_5=21|346$
are not. Thus, $\{S_1,\ldots, S_5\}$
is not  weakly compatible.
Note that, like in the case of (full) splits, it is easy to
see that any collection of pairwise compatible  partial
splits is also  weakly compatible.

Clearly any three partial splits $S_i=A_i|\tilde{A_i}\in\s(X)$, $i=1,2,3$,
for which precisely one of the four intersections in
(\ref{thinly}) is empty also satisfies Condition (\ref{SK-rule}). 
Thus $\theta_Y$
may be applied to $S_1,S_2,S_3$. However, as
the example of the
set $\{127|3456, 1234|567, 235|146\}$ shows,
application  of $\theta_Y$ to a
 weakly compatible collection of partial splits does not, in
general, yield a  weakly compatible collection of partial splits.
Also it should be noted that
$\theta_M$ applied to a  weakly compatible collection of partial splits
 does not always
yield a  weakly compatible collection of partial splits.

However, the next result whose proof is straight forward holds.

\begin{lemma}\label{pink}
Suppose $\s,\s'\subseteq\s(X)$. If $\s'$ is
 weakly compatible and $\s\preceq \s'$, then $\s$ must also be
 weakly compatible.
\end{lemma}

\subsection{Circular collections of partial splits}\label{circular-splits}

We now turn our attention to the study of
a special class of
weakly compatible collections of partial splits called
circular collections of partial splits. To be able to
state their definition, we require some more terminology which
we introduce next.

A {\em cycle} $C$ is a connected graph with $|V(C)|\geq 3$ and every
vertex has degree 2. We call $C$ an {\em $X$-cycle} 
if the vertex set of $C$ is $X$.
For $x_i\in X$ ($1\leq i\leq n:=|X|$) and $C$ an $X$-cycle,
we call $x_1,x_2,\ldots,x_n, x_{n+1}=x_1$ a
{\em vertex ordering (of $C$)} if the edge set of $C$ coincides with the
set of all $2$-sets $\{x_i,x_{i+1}\}$ of $X$, $i=1,\ldots,n$.

For a graph $G=(V,E)$ and some subset $E'$ of $E$, we denote by $G-E'$ the
graph obtained from $G$ by deleting the edges in $E'$.
We say that a partial $X$-split $A|\tilde{A}$ is {\em displayed} by an
$X$-cycle $C$ if there exist two distinct edges $e_1$ and $e_2$ in
$C$ such that the vertex set of
one of the two components of $C-\{e_1,e_2\}$ contains
$A$ and the other one contains $\tilde{A}$. More generally,
we say that a set $\s\subseteq \s(X)$ of partial
splits is {\em displayed} by an $X$-cycle $C$ if every partial split in $\s$ is
displayed by $C$. Finally, we say that a collection $\s\subseteq \s(X)$
is {\em circular} if there exists some $X$-cycle $C$ such that
every partial split in $\s$ is displayed by $C$. Note that every
split collection
in $\s(X)$ displayed by a circular phylogenetic network is circular.

As is well-known, every circular split system is in
particular weakly compatible. The next result shows that an analogous
result holds for collections of partial splits.

\begin{lemma}\label{circular-twc}
Suppose $\s\subseteq \s(X)$. If $\s$ is circular then $\s$ is also
 weakly compatible.
\end{lemma}
\pf
Suppose $C$ is an $X$-cycle that displays $\s$
but there exist three partial
splits $S_1,S_2,S_3\in \s$ such that
with $A_i\in S_i$, $i=1,2,3$,
playing the role of their namesakes in (\ref{thinly})
none of the four
intersections in (\ref{thinly}) is empty.
Then $S_1$ and $S_2$ are incompatible and, since $S_1$ and $S_2$
are displayed by $C$,
there must exist edges $e_1,e_1',e_2,e_2'\in E(C)$
such that, for all $i,j\in \{1,2\}$
distinct, the vertex set of one component of $C-\{e_i,e_i'\}$ contains
$A_i \cup e_j$ and the other
contains $\tilde{A_i} \cup e_j'$.
Since $S_3$ is displayed by $C$  and neither
$A_1\cap A_2\cap A_3$ nor $\tilde{A_1}\cap \tilde{A_2}\cap A_3$,
nor $\tilde{A_1}\cap A_2\cap \tilde{A_3}$ is empty, it follows that
$A_1\cap \tilde{A_2}\cap\tilde{ A_3}=\emptyset$, which is impossible.
\epf

As in the case of full splits, the converse
of the above lemma is not true in general. For example, the
set $\s$ comprising the partial splits $S_1=12|35$,
$S_2=125|34$,
$S_3=13|245$ and $S_4=135|24$ is  weakly compatible since
the sets $\{S_1,S_2\} $ and
$\{S_3,S_4\} $ are pairwise compatible. Yet, as can be easily checked,
 $\s$ is not circular.

Corresponding to Lemma \ref{pink}, we have:

\begin{lemma}\label{purple}
Suppose $\s,\s'\subseteq\s(X)$. If $\s'$ is
displayed by an $X$-cycle $C$ and $\s\preceq \s'$, then $\s$ is
also displayed by $C$.
\end{lemma}

\subsection{Circularity and the $M$- and $Y$-rule}

As was noted earlier, neither the $Y$-rule nor the $M$-rule
preserve weak compatibility in general. As the next result shows,
 the situation is different for the special case of
circular collections of partial splits.


\begin{proposition}\label{circ}
Suppose  $\s,\s'\in \p(X)$
and $C$ is an $X$-cycle.
If $\s'$ is obtained from $\s$ by a single application of either
$\theta_Y$ or $\theta_M$
then $\s$ is displayed by $C$ if and only if $\s'$ is displayed by $C$.
\end{proposition}
\pf Suppose $\s,\s'\in \p(X)$ and $C$ is an $X$-cycle. We start the proof
with noting
that, regardless of whether $\s'$ is obtained from a single
application of either $\theta_Y$ or $\theta_M$
to $\s$, $\s$ is displayed by $C$ whenever $\s'$ is
displayed by $C$ in view of Lemma~\ref{purple}.

Conversely, suppose that $\s$ is displayed by $C$.
Assume first that
 $\s'$ is obtained from $\s$ by a single application of $\theta_Y$. Let
$\{S_1,S_2,S_3\}\subseteq \s$ be the set to which  $\theta_Y$
is applied. With $A_i\in S_i$, $i=1,2,3$,
playing the role of their namesakes in the statement of (\ref{SK-rule}),
we may assume without loss
of generality that none of the three intersections
$D_1=A_1\cap A_2\cap A_3$,  $D_2=\tilde{A_1}\cap \tilde{A_2}\cap A_3$, and
 $D_3=\tilde{A_1}\cap A_2\cap \tilde{A_3}$ is empty but that
$A_1\cap \tilde{A_2}\cap\tilde{A_3}=\emptyset$.  It
suffices to show that the partial split
$S= A_3\cup (A_1\cap \tilde{A_2})|\tilde {A_3}$ is displayed by $C$.

\begin{figure}[h]
\begin{center}
\input{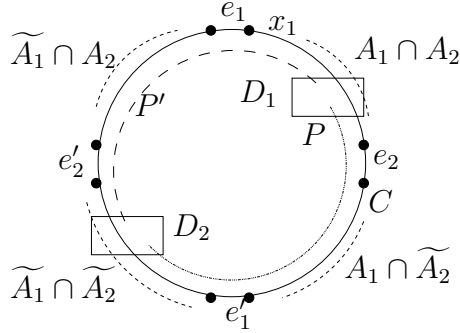}
\caption{\label{proof}
A schematic representation of the two alternative locations
for $\tilde{A_3}$ (cf proof of Proposition~\ref{circ}). The
closed curve is the $X$-cycle $C$,
the four curves with the short dashes
represent the four non-empty intersections
$A_1\cap A_2$, $A_1\cap \tilde{A_2}$, $\tilde{A_1}\cap \tilde{A_2}$
and $\tilde{A_1}\cap A_2$ (note that each of them can consist
of more than one part), the rectangles mark the intersections $D_1$ and $D_2$,
and the dotted and dashed curves
represent the two paths $P$ and $P'$ on $C$ on which
$\tilde{A_3}$ can lie.
}
\end{center}
\end{figure}

Clearly, if $A_1\cap \tilde{A_2}=\emptyset$ then $S=S_3$
and, therefore, $S$ is displayed by $C$.
So assume $A_1\cap \tilde{A_2}\not=\emptyset$.
Then since by assumption $D_i\not=\emptyset$, $i=1,2,3$,
and  $S_1$
and $S_2$ are displayed by  $C$, there must exist four distinct
edges $e_1,e_1',e_2,e_2'\in E(C)$ such that, for all $i,j\in \{1,2\}$
distinct, one component of $C-\{e_i,e_i'\}$ contains
$A_i$ in its vertex set and $e_j\subseteq A_i$ and the other
contains $\tilde{A_i}$ in its vertex set and $e_j'\subseteq \tilde{A_i}$.
Without loss of generality, we may assume that
$X=\{x_1,\ldots, x_n\}$, $n\geq 3$, that $x_1,x_2,\ldots, x_n$ is a
vertex ordering of $C$, and that $e_1=\{x_n,x_1\}$. Furthermore, we
may also assume without loss of generality that
the component of  $C-\{e_1,e_1'\}$ that contains $x_1$ in its
vertex set also contains
$A_1$. Since $D_1\not=\emptyset\not=D_2$,
and $S_3$ is displayed by
$C$ there must exist distinct paths $P$ and $P'$ in $C$ such that
either $\tilde{A_3}\subseteq V(P) $
or $\tilde{A_3}\subseteq V(P')$ (see Figure~\ref{proof}).
If $\tilde{A_3}\subseteq V(P)$ then
$$
\emptyset = \tilde{A_1}\cap A_2\cap \tilde{A_3}=D_3
$$
which is impossible.
Thus $\tilde{A_3}\subseteq V(P')$ must hold. Suppose $y,z\in V(C)$ are
such that when starting at $x_1$ and traversing $C$ clockwise $y$
is contained in $A_3$ and the next vertex $y'$ on $C$ with
$y'\in A_3\cup\tilde{A_3}$ is contained in  $\tilde{A_3}$
whereas $z\in \tilde{A_3}$ and the next vertex $z'$ on $C$
with $z'\in A_3\cup\tilde{A_3}$ is contained in  $A_3$.

 Let $P''$ denote the path from $z'$ to $y$ (taken clockwise). Then $e_2$ and
$e_1'$ are edges on $P''$ and so $A_1\cap \tilde{A_2}\subseteq V(P'')$.
The choice of $y$ and $z'$ implies
$V(P'')\cap\tilde{A_3}=\emptyset$ and
$A_3\cup (A_1\cap \tilde{A_2}) \subseteq V(P'')$.
Hence, the split $V(P'')|X-V(P'')$ which is displayed
by $C$ extends the partial split $S$.
Thus $C$ displays $S$. This concludes the proof in case
the applied closure rule applied is $\theta_Y$.

To conclude the proof of the proposition suppose
$\s'$ is obtained from $\s$ by a single application of $\theta_M$.
Let  $\{S_1,S_2\}\subseteq \s$ be the set to which $\theta_M$
is applied. With $A_i\in S_i$, $i=1,2$, we may assume
without loss of generality that
$A_1\cap A_2\not=\emptyset$ and
$\tilde{A_1}\cap \tilde{A_2}\not=\emptyset $.
If $\theta_M$ applies trivially to $\s$ then $\s=\s'$ and so
$\s'$ must be displayed by $C$. If $\theta_M$ does not apply trivially to
$\s$ it suffices to show
 that  $C$ displays $(A_1\cap A_2)|(\tilde{A_1}\cup \tilde{A_2})$.

Since $S_1$ and $S_2$ are displayed by $C$ there must exist edges
$e_i,e_i'\in E(C)$ such that the vertex set of one of the two components
$P_i,P_i'$ of $C-\{e_i,e_i'\}$ contains $A_i$ and the other
 contains $\tilde{A_i}$, $i=1,2$.
Put $k:=| \{e_1,e_1'\}\cap \{e_2,e_2'\} |$ and note that $0\leq k\leq 2$.
Without loss of generality, we may
assume $A_i\subseteq V(P_i)$ and $\tilde{A_i}\subseteq V(P_i')$, $i=1,2$.
Then, $\emptyset\not= A_1\cap A_2\subseteq V(P_1)\cap V(P_2)$.
Since $V(P_1)\cap V(P_2)$ is the vertex set of one of the $4-k$
components of $C$ with the edges $e_i,e_i'$, $i=1,2$
removed, it follows that there must exist
two distinct edges $e_3,e_4$ among the edges $e_1,e_1',e_2,e_2'$
so that the vertex set of one of the two components of $C-\{e_3,e_4\}$
is $V(P_1)\cap V(P_2)$. Since
$$
X-(V(P_1)\cap V(P_2))=
(X-V(P_1))\cup (X-V(P_2))=V(P_1')\cup V(P_2')
$$
is the vertex set of the other component of
$C-\{e_3,e_4\}$
and $\tilde{A_i}\subseteq V(P_i')$, $i=1,2$, it follows that $C$ displays
$(A_1\cap A_2)|(\tilde{A_1}\cup \tilde{A_2})$. This concludes the
proof in case $\s'$ is obtained from $\s$ by a
single application of $\theta_M$ and thus the proof of the proposition.
\epf

Interestingly, the $Z$-rule
does not preserve circularity in general. An
example in point is the $X$-cycle $C$ with
$X=\{1,\ldots, 5\}$ and the natural ordering of the
elements of $X$ as vertex ordering. Then the partial splits $S_1=13|45$
and $S_2=34|25$ are clearly displayed by $C$. Yet
the $Z$-rule applied to $\{S_1,S_2\}$
generates the partial splits $13|245$ and $25|134$
which cannot be displayed by $C$.

\section{Split closure sequences and split closures}\label{spl-sequence}

In this section, we associate to a set $\s$ of
partial splits a split closure sequence and define the last element
of such a sequence to be a split closure of $\s$. We also establish
a key result for this paper which shows that under certain
circumstances a split closure is unique.

\subsection{Split closure sequences}\label{scs}

Suppose $\s\in \p(X)$ is a collection of partial splits
that satisfies some partial splits property $(P)$
such as, for example,  weak compatibility
and $\theta$ is one of the closure rules considered in this paper. Following
\cite{SS01}, we associate a split closure sequence $\sigma$
and a split closure to $\s$ as follows.
$$
\sigma: \s_0,\s_1,\s_2,\ldots,\s_i,\s_{i+1}, \ldots
$$
is a strictly increasing (with respect to $\preceq$)
sequence of sets in $\cP(X)$ so that
$\s=\s_0$ and, for all $i\geq 1$, $\s_{i+1}$ is obtained by one
non-trivial application of $\theta$ to $\s_i$
whenever $\s_i$ satisfies $(P)$. Note that
since $X$ is finite,
there must exist a last element $\s_n$ in $\sigma$
such that $\s_n$ either satisfies $(P)$ and is closed
under $\theta$ or
$\s_n$ does not satisfy $(P)$. In the latter case we reset
$\s_n$ to be a new element
$\omega \notin \cP(X)$.
We refer to $\sigma$ as a {\em split closure sequence}
for $\s$ and call $n$ the length of
$\sigma$. In addition, we call the last element
of $\sigma$ a  {\em split closure} of $\s$. Note that
in case $\s_n\not=\omega$,
$\theta$ applies only trivially to $\s_n$.

The following combinations of $(P)$ and $\theta$ are
of interest to us: 
\begin{enumerate}
\item[(a)] $(P)$ is the property that $\s$ is
weakly compatible  and $\theta$ is the $Y$-rule.
\item[(b)] $(P)$ is unspecified and
$\theta$ is the $M$-rule.
\item[(c)]  $(P)$ is the property that $\s$ is
weakly compatible and $\theta$ is the $M/Y$-combination closure rule
$\theta_{M / Y}$
which applies $\theta_M$ or $\theta_Y$ to $\s$.
\end{enumerate}

To elucidate the notion of a split closure sequence and
a split closure associated to a set in $\p(X)$ we next present an
example for the assignments of $(P)$ and $\theta$ specified in (a).
Consider the set   $X=\{1,2,3,4,5\}$ together with
the collection $\s$ comprising of the partial $X$-splits
$S_1=12|34$, $S_2=23|14$, $S_3=15|24$, and $S_4= 45|13$.
Clearly, $\s$ is displayed by an $X$-cycle $C$ with
vertex ordering $1,2,3,4,5$. Thus $\s$ is circular and so, by
Lemma~\ref{circular-twc}, $\s$
is weakly compatible. Now $\theta_Y$ applied to $\{S_1,S_2,S_3\}$
generates the
split $S_3'=15|234$, $\theta_Y$ applied to $\{S_1,S_2,S_4\}$ generates the
split $S_4'=45|123$ and $\theta_Y$ applied to $\{S_2,S_3',S_4'\}$ generates the
split $S_2'=145|23 $. Since every subset of
$\s'= \{S_1,S'_2,S'_3,S_4'\}$ of size three contains two pairwise
compatible full splits, $\theta_Y$ can only be applied trivially to $\s'$.
Hence, the sequence $S_0=\s$, $\s_1=\{S_1,S_2,S'_3,S_4\}$,
$\s_2=\{S_1,S_2,S'_3,S'_4\}$,
$\s'$  is a split closure sequence for $\s$
of length $3$
and $\s'$ is a split closure for $\s$.

Regarding (c), it should be noted that even if for some
$\s\in \p(X)$ two distinct split closure sequences have
the same length and terminate in the same element $\s'\not=\omega$
one of them might
utilise fewer applications of $\theta_Y$ (and thus
more applications of $\theta_M$!)
than the other. For the previous example, one
way to construct two such sequences is to
exploit the following relationship between
the $Y$-rule and the $M$-rule for
$\{S_2, S_3', S_4'\}$.


\begin{proposition}\label{relationship}
Suppose $\s=\{S_i=A_i|\tilde{A_i}:i=1,2,3\}\in \p(X)$ is 
such that
$A_1\subseteq A_2$ and
$\tilde{A_2}-\tilde{A_1}\subseteq \tilde{A_3}
\subseteq\tilde{A_1}\cup \tilde{A_2}$.
If  the $Y$-rule applies
to $\s$  then
$$
\theta_{Y}(\s)^-=\{\tilde{A_1}\cup \tilde{A_2}|A_1, S_2, S_3\}
=\{S_3\}\cup \theta_{M}(S_1,S_2)^-.
$$
\end{proposition}
\pf Assume that $\s$ and  $S_i$ and $A_i$, $i=1,2,3$, are such that
the assumptions of the proposition are satisfied. Then
$A_1\cap \tilde{A_2}=\emptyset$. Combined with the assumption that 
$\theta_Y$ applies to
$\s$, it follows that either (\ref{SK-rule}) is satisfied
 with $A_i$, $i=1,2,3$ playing
the roles of their namesakes in the statement of
(\ref{SK-rule}) or (\ref{SK-rule}) is
satisfied with $A_3$ playing the role of $\tilde{A_3}$ and
$A_i$ playing the role of $A_i$, $i=1,2$, in that statement.
But the latter
alternative cannot hold since this implies
$A_1\cap A_2\cap \tilde{A_3}\not=\emptyset$ whereas
the assumption $\tilde{A_3}\subseteq\tilde{A_1}\cup \tilde{A_2}$ implies
$$
A_1\cap A_2\cap \tilde{A_3}
\subseteq A_1\cap A_2\cap (\tilde{A_1}\cup \tilde{A_2})=
(A_1\cap A_2\cap \tilde{A_1})\cup  (A_1\cap A_2\cap \tilde{A_2})
=\emptyset.
$$
Hence, (\ref{SK-rule}) is satisfied with $A_i$, $i=1,2,3$ playing
the roles of their namesakes in the statement of (\ref{SK-rule}).
Let $S_i'$, $i=1,2,3$ be as in
the statement of the $Y$-rule. Then
$\tilde{A_2}-\tilde{A_1} \subseteq \tilde{A_3} $ implies
$S_1'=\tilde{A_1}\cup \tilde{A_2}|A_1$. And since
$A_1\subseteq A_2$, we have $S_2'=S_2$ and $A_1\cap \tilde{A_2}=\emptyset$
which in turn implies $S_3'=S_3$. Consequently,
$$
\theta_{Y}(\s)^-=\{\tilde{A_1}\cup \tilde{A_2}|A_1, S_2, S_3\}.
$$

To observe the remaining set equality, note
that since (\ref{SK-rule}) is satisfied with $A_i$, $i=1,2,3$ playing
the roles of their namesakes in the statement of (\ref{SK-rule})
neither
$A_1\cap A_2$ nor $\tilde{A_1}\cap \tilde{A_2}$ can be empty.
Let $S_1'$ be as in the statement of the $M$-rule. Then
$A_1\subseteq A_2$
implies $S_1'=\tilde{A_1}\cup \tilde{A_2}|A_1$ and
$S_2'=\tilde{A_1}\cap \tilde{A_2}|A_2$.
This implies the sought after set equality and thus proves the proposition.
\epf

Clearly independent of which one of the rules $\theta_Y$, $\theta_M$
or $\theta_{M/Y}$ is applied, a split closure sequence must always be finite
since $X$ is finite. In addition and by
applying the same arguments as Semple and Steel in \cite{SS01} one can
show that, for the assignments of $(P)$ and $\theta$ as described in (a),
the length of a split closure sequence for a  weakly compatible
set $\s\in \cP(X)$ is bounded from above by
$|\s|\cdot|X| -\Sigma_{\{A,B\}\in\s}|A\cup B| $.

\subsection{Split closures}

We start with a lemma that is crucial for
showing that the split closure of some collection $\s\in \p(X)$
is unique in any of the three combinations for ($P$)
and $\theta$ stated in (a) -- (c).

\begin{lemma}\label{key}
Suppose $\s\in \p(X)$, $\overline{\s}\not=\omega$
is a split closure of $\s$ and $\s_r$ and $\s_{r+1}$
are two consecutive elements in a split closure sequence for $\s$.
\begin{enumerate}
\item[(i)] If $\s_r$ is  weakly compatible,
$ \s_r\preceq \overline{\s}$ and  $\s_{r+1}$ is obtained from
$\s_r$ by one application of $\theta_Y$, then
$\s_{r+1}$ is  weakly compatible
and $ \s_{r+1}\preceq \overline{\s}$.
\item[(ii)]  If $\s_{r+1}$ is obtained from
$\s_r$ by one application of $\theta_M$
and $\s_r\preceq \overline{\s}$,
then  $ \s_{r+1}\preceq \overline{\s}$.
\end{enumerate}
\end{lemma}
\pf Suppose $\s,\overline{\s},\s_r, \s_{r+1}$ are as
in the statement of the lemma.

(i) Assume $\s_{r+1}$ is obtained
from $\s_r$ by applying $\theta_Y$ to some set $\{S_1,S_2,S_3\}$
contained in $ \s_r$.
For $i=1,2,3$ and with $A_i\in  S_i$ playing the role of
their namesakes in the statement of (\ref{SK-rule}), we obtain
$$
S_3'=A_3\cup (A_1\cap \tilde{A_2})|\tilde{A_3},
S_2'=A_2'\cup(A_1\cap \tilde{A_3})|\tilde{A_2}\mbox{, and }
S_1'=\tilde{A_1}\cup(\tilde{A_3}\cap \tilde{A_2})|A_1.
$$
It follows that
$$
\s_{r+1}=(\s_r\cup\{S_1', S_2',S_3'\})^-.
$$
Since $\s_r\preceq \overline{\s}$,
there exist partial splits
$S_i''=A_i''|\tilde{A_i''}\in \overline{\s}$ with $S_i''$ extending $S_i$,
$i=1,2,3$. Without loss of generality we may assume for all $i$ that
$A_i\subseteq A_i''$ and $\tilde{A_i}\subseteq \tilde{A_i''}$.
Since $\overline{\s}$ is  weakly compatible, (\ref{SK-rule}) 
is satisfied by
$\{S_1'', S_2'', S_3''\}$ with
$A_1''\cap \tilde{A_2''}\cap\tilde{A_3''} =\emptyset $.
Since,  by assumption, $\overline{\s}\not =\omega$ and so $\theta_Y$
 applies trivially to $\overline{\s}$, we must have
$A_1''\cap \tilde{A_2''}\subseteq A_3''$,
$A_1''\cap \tilde{A_3''}\subseteq  A_2''$,
and $\tilde{A_3''}\cap \tilde{A_2''}\subseteq  \tilde{A_1''}$. It follows
that,
for all $i=1,2,3$, $S_i'$ is extended by $S_i''$
which in turn implies  $\s_{r+1}\preceq \overline{\s}$.
Since $\overline{\s}\not=\omega$ and so
$\overline{\s}$ is  weakly compatible, Lemma~\ref{pink}
implies that $\s_{r+1}$ is  weakly compatible.

(ii): Suppose
$\s_{r+1}$ is obtained from
$\s_r$ by applying $\theta_M$ to some set $\{S_1,S_2\}\subseteq \s_r$.
Put $S_i=A_i| \tilde{A_i}$,
$i=1,2$, and assume without loss of  generality, that
$A_1\cap A_2\not=\emptyset$ and $\tilde{A_1}\cap \tilde{A_2}\not=\emptyset$.
Then
$$
S_1'=A_1\cap A_2|\tilde{A_1}\cup\tilde{A_2''}\mbox{ and }
S_2'=A_1\cup A_2|\tilde{A_1}\cap\tilde{A_2''}
$$
and so
$$
\s_{r+1}=(\s_r\cup\{S_1',S_2'\})^-.
$$
By assumption,  $\s_r\preceq \overline{\s}$
and so there exist  partial splits
$S_1''=A_1''|\tilde{A_1''}$, $S_2''=A_2''|\tilde{A_2''}$
in $\overline{\s}$ with $S_i''$ extending $S_i$,
$i=1,2$. Without loss of generality we may assume
$A_i\subseteq A_i''$ and $\tilde{A_i}\subseteq \tilde{A_i''}$,
$i=1,2$.  Then
$A_1''\cap A_2''\not=\emptyset$ and
$\tilde{A_1''}\cap \tilde{A_2''}\not=\emptyset$.
Since, by assumption,
$\overline{\s}\not=\omega$ and so $\theta_M$ only applies trivially to
 $\overline{\s}$ we have $\theta_M(S_1'',S_2'')\preceq \overline{\s}$.
Hence, there exist partial splits $S_j=A_j|\tilde{A_j}\in \overline{\s}$,
$j=3,4$, so that $A_1''\cap A_2''|\tilde{A_1''}\cup \tilde{A_2''}$
is extended by $S_3$ and
$\tilde{A_1''}\cap \tilde{A_2''}|A_1''\cup A_2''$ is extended by $S_4$.
Without loss of generality, we may assume that
$A_1''\cap A_2''\subseteq A_3$ and
$\tilde{A_1''}\cup \tilde{A_2''}\subseteq \tilde{A_3}$
and that $A_1''\cup A_2''\subseteq A_4$ and
$\tilde{A_1''}\cap \tilde{A_2''}\subseteq \tilde{A_4}$.
Then
$$
A_1\cap A_2\subseteq A_1''\cap A_2''\subseteq A_3
\mbox{ and } \tilde{A_1}\cup
\tilde{A_2}\subseteq \tilde{A_1''}\cup \tilde{A_2''}\subseteq \tilde{A_3''}
$$
and so  $S_3$ extends $S_1'$. Similarly, it follows that
$S_4$ extends $S_2'$. Thus,
 $\s_{r+1}\preceq \overline{\s}$.
\epf

With this result in hand, we are now in the position to present a key result.

\begin{theorem}\label{split-closure}
Suppose $\s\in \p(X)$.
Then any two split closures
for $\s$ are the same if
\begin{enumerate}
\item[(i)] $\s$ is  weakly compatible and solely
the $Y$-rule is used to obtain a split closure for $\s$,
\item[(ii)] solely the
$M$-rule is used to obtain a split closure for $\s$, or
\item[(iii)] $\s$ is  weakly compatible
and solely the $M/Y$-rule
is used to obtain a split closure for $\s$.
\end{enumerate}
\end{theorem}
\pf Suppose $\s\in \p(X)$.
We start with remarking that we prove Statements (i), (ii)
and (iii) collectively as the proof of all three statements
relies on an inductive argument on the length of a split closure sequence
for $\s$. However, since the arguments for the inductive step
differ under the assumptions made in (i), (ii) and (iii), we
discuss each inductive step separately.

Suppose that the assumptions made in
(i) or in (ii) or in (iii) hold.
If every split closure of $\s$ is $\omega$ then the
theorem holds trivially. So we may assume that there exists a
split closure $\overline{\s}$ of $\s$ with $\overline{\s}\not=\omega$.
We proceed by showing that every other split closure of $\s$ must equal
$\overline{\s}$. Suppose that $\sigma: \s_0=\s,\s_1,\s_2,\ldots,\s_n$
is a split closure sequence of $\s$. We now
use induction on $n$ to show that if $\s$ satisfies the assumptions made:
\begin{enumerate}
\item[] in (i) then, for all $i\in\{0,1,\ldots, n\}$,
\begin{eqnarray}\label{closure-claim1}
\s_i\mbox{ is  weakly compatible and } \s_i\preceq \overline{\s};
\end{eqnarray}
\item[] in (ii) then, for all $i\in\{0,1,\ldots, n\}$,
\begin{eqnarray}\label{closure-claim2}
\s_i\preceq \overline{\s};
\end{eqnarray}
\item[] in (iii) then, for all $i\in\{0,1,\ldots, n\}$,
\begin{eqnarray}\label{closure-claim3}
\s_i\mbox{ is  weakly compatible}
 \mbox{ and } \s_i\preceq \overline{\s}.
\end{eqnarray}
\end{enumerate}

We start with assuming that $\s$ satisfies the assumptions made
in (i), that is,
$\s$ is  weakly compatible and solely the $Y$-rule
is used to generate the elements of $\sigma$.  If $i=0$  then
(\ref{closure-claim1}) obviously holds since then $\s_i=\s_0$
and $\s_0$ satisfies the properties stated in (\ref{closure-claim1}).
Now suppose that (\ref{closure-claim1})
holds for some $i\in\{0,1,\ldots,n-1\}$. Then,
by Lemma~\ref{key}(i),
$\s_{i+1}$ is  weakly compatible
and $ \s_{i+1}\preceq \overline{\s}$.
This completes the induction step and thereby establishes
(\ref{closure-claim1}).

Next, assume that only the $M$-rule is used
to generate the elements in
$\s$. If $i=0$, then (\ref{closure-claim2}) holds  since
then $\s_i=\s_0$ and $\s_0$ satisfies (\ref{closure-claim2}). Assume that
(\ref{closure-claim2})
holds for some $i\in\{0,1,\ldots,n-1\}$.
Then Lemma \ref{key}(ii) implies $\s_{r+1}\preceq
\overline{\s}$  which completes the induction step and thereby
establishes (\ref{closure-claim2}).

Finally, assume that $\s$ satisfies the
assumptions made in (iii), that is, $\s_r$ is
weakly compatible and only the $M/Y$-rule
is used to generate the elements in $\sigma$.
If $i=0$,  then
(\ref{closure-claim3}) obviously holds since then $\s_i=\s_0$
and $\s_0$ is closed under 
$\theta_{M/Y}$.
Now suppose that (\ref{closure-claim3})
holds for some $i\in\{0,1,\ldots,n-1\}$. Then
$\s_i$ is  weakly compatible, $\s_r\preceq \overline{\s}$, and
one of the
following two cases must hold. Either
(a) $\s_{r+1}$ is obtained from $\s_r$ by applying $\theta_Y$ or 
(b) $\s_{r+1}$ is
obtained from $\s_r$ by applying $\theta_M$.

If Case (a) holds, the proof of the inductive step in (i)
implies that $\s_{i+1}$ is  weakly compatible and
$\s_{i+1}\preceq \overline{\s}$.

If Case (b) holds, $\s_{i+1}\preceq \overline{\s}$
follows from the proof of the inductive step in (ii). That
$\s_{i+1}$ is  weakly compatible follows from
Lemma~\ref{key}(ii) and the fact that
$\overline{\s}$ is weakly compatible.
This completes the induction step and thereby
establishes~(\ref{closure-claim3}).

We conclude with noting that for $i=n$, we obtain $\s_n\preceq \overline{\s}$
regardless of whether we are assuming (i) or (ii) or (iii) to hold.
In case of (i)  holding and
applying (\ref{closure-claim1}) to $i=n$
or (iii)  holding and applying (\ref{closure-claim3}) to $i=n$,
we see that $\s_n$ is
 weakly compatible. By interchanging the roles of
$\s_n$ and $\overline{\s}$, we deduce $\overline{\s}\preceq \s_n$.
Thus, under the assumptions made in (i) or (ii) or (iii), we have
$\overline{\s}= \s_n$ which concludes the proof of the theorem.
\epf

Extending in the case of ($P$) denoting the
condition ``$\s$ is  weakly compatible'' and $\theta$ denoting
either the $Y$-rule or the $M/Y$-rule,
the definition of the split closure to non  weakly compatible
sets in $ \cP(X)$
by defining the split closure of such sets to be $\omega$, we obtain

\begin{corollary}\label{summary}
Suppose $\theta$ is either the $Y$- or $M$- or $M/Y$-rule.
Then for any $\s\in \p(X)$,
any two split closures for $\s$
obtained via $\theta$ are the same.
\end{corollary}

Bearing in mind Corollary~\ref{summary}, we denote
for
$\theta\in \{\theta_Y,\theta_M, \theta_{M/Y}\}$
the split closure of
a set $\s\in \p(X)$ by
$\langle\s\rangle_{\theta}$, Note that
 $|\langle \s\rangle_{\theta_Y}|\leq |\s|$
but that neither $|\langle \s\rangle_{\theta_M}|\leq |\s|$
nor $|\langle \s\rangle_{\theta_{M/Y}}|\leq |\s|$ have to hold.
Also note that, if we  denote the collection
$\p(X)\cup\{\omega\}$ by $\p_{\omega}(X)$, define
$\langle\omega \rangle_{\theta}=\omega$ for some closure rule $\theta$,
and put $\s\preceq\omega$ for all $\s\in \cP(X)$,
then  the
split closure with respect to $\theta\in
\{\theta_Y,\theta_M,M / Y\}$ satisfies
the usual properties of a closure operation. More precisely,
for all $\s,\s'\in \p_{\omega}(X)$ we have
$\s\preceq \langle \s \rangle_{\theta}$, if $\s\preceq\s'$ then
$\langle \s \rangle_{\theta} \preceq  \langle \s'\rangle_{\theta} $, and
$ \langle\langle \s \rangle_{\theta} \rangle_{\theta}
=\langle \s\rangle_{\theta}$.

As an immediate consequence of
Lemma~\ref{circular-twc}, Proposition \ref{circ}, and
Theorem \ref{split-closure}, we obtain our main result which
we state next.

\begin{corollary} \label{cyclosure}
Suppose $\s\in \p(X)$ and $C$ is an $X$-cycle. Then $\s$ is
displayed by $C$ if and only if  $\langle \s\rangle_{\theta_{M/Y}}$
is displayed by $C$. In that case $\langle \s\rangle_{\theta_Y}$ and
$\langle \s\rangle_{\theta_M}$ are also displayed by $C$.
\end{corollary}

We conclude this section with remarking that, in general,
 not all elements
in $\langle\s\rangle_{\theta_{M/Y}}$ need to be full splits on $X$,
$\s\in \cP(X)$ circular. However, it is reasonable to assume that
those that have been extended to full splits on $X$ contain
phylogenetically relevant information and programs such as
e.\,g.\,SplitsTree4 \cite{HB05}
may be employed to produce a circular phylogenetic network
that displays them.
For the following we refer to the combination of the $M/Y$-rule with
a phylogenetic network generation package such as SplitsTree4
as the {\em $MY$-closure approach}. Although a detailed 
analysis of this approach is beyond the scope of this
paper and will be presented elsewhere, we note that 
the MY-closure approach cannot be polynomial in the worst case since
if the collection of partial splits  comprises of all 
$ 3{n \choose 4}$ partial
$X$-splits $A|B$ with $|A|=2=|B|$ and  $n=|X|$ 
then $\theta_{M/Y}$  will generate all $2^{(n-1)}-1$ splits of $X$



\section{An example: The ring of life}\label{rivera-exmple}

\begin{figure}[h]
\begin{center}
\input{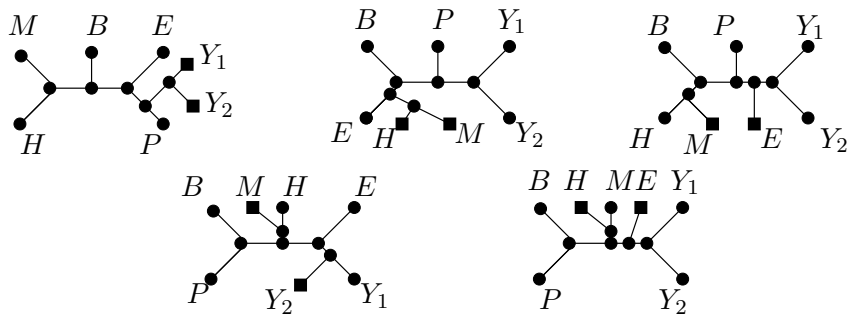}
\caption{\label{ring}
Five most probable phylogenetic trees that appeared in slightly
different form in \cite[Fig. 1]{RL04} (cf text for details).
}
\end{center}
\end{figure}

One of the most fiercely debated questions amongst biologists is the origin
of eukaryotes (essentially cells that have a nucleus and organelles)
\cite{ME04}. The main reason for this is that
eukaryotes have eubakteria-like genes as well as archaebacteria-like genes
making it very difficult to establish the evolutionary relationships
between eukaryotes and prokaryotes (essentially cells that lack
nucleus and organelles) which is the collective name for
eubakteria and archaebacteria.
To help shed light into this question, Rivera et al.
\cite{RL04} analysed 10 bacterial genomes. The 5 most probable
phylogenetic trees resulting from their analysis are presented in Fig.~1
of that paper. For the convenience of the reader, we
depict them in slightly different form in Fig.~\ref{ring}. Note that
the collection of splits displayed by these trees is circular
and also that, when
ignoring the fact that the leaves are marked with
different symbols, the last 2 trees are the same.

\begin{figure}[h]
\begin{center}
\input{MY-net.pstex_t}
\hfill
\raisebox{-8ex}
{\includegraphics[scale=0.36
]{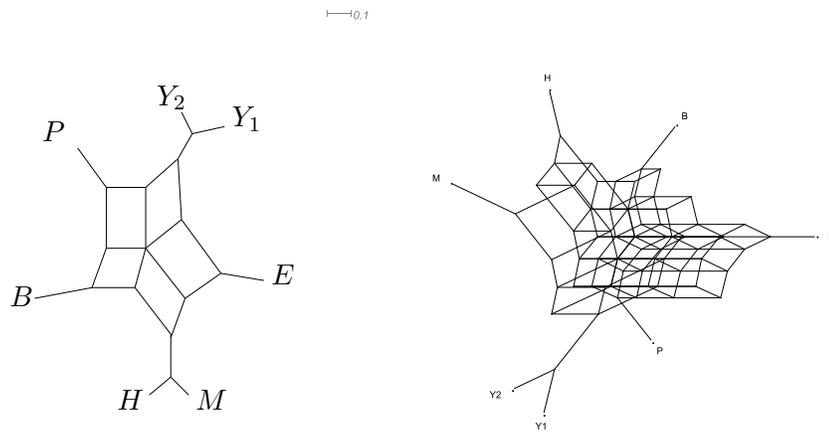}}
\caption{\label{lake-network}
Left, a circular network on 2 yeast genomes,
an $\alpha$-probacterium,
a bacillus,
a halobacterium,
a methnaococcus,
an ecocyte,
and an  archaeoglobium (the genome abbreviations
follow \cite{RL04}). It displays
the split collection
inferred from the collection of partial splits induced by the trees in
Fig.~\ref{ring} with the leaves marked with a square
plus their incident edges removed and the resulting degree
2 vertices suppressed using the $MY$-closure approach. Right, the Z-closure
super-network on the same set of partial splits}
\end{center}
\end{figure}

Using a technique called {\em Conditioned
Reconstruction} \cite{RL04}, Rivera et al. constructed the phylogenetic
network depicted in
Fig.~\ref{lake-network} with the degree 5 interior vertex plus all its
incident edges removed and all resulting degree 2 vertices
suppressed. The resulting structure they then interpreted as
lending support to the idea that, in its early stages, evolution
was not tree-like but rather more like a ring (hence the 
term ``ring of life'') with the
eukaryotic genome being the result of a fusion of 2 diverse
procaryotic genomes~\cite{RL04}.


To find out how dependent Rivera et al.'s ring
of life is on the fact that all 5 trees are on the same leaf set, we
randomly removed pairs of leaves plus their incident
edges (suppressing resulting degree 2 vertices and
always ensuring that there were no 2 trees from which the same pair of leaves
was removed) resulting in 5 trees $T_1,\ldots,T_5$
on 5 leaves.
Perhaps not surprisingly, we found that, in general, removal
of pairs of leaves did not allow us to recover Rivera et al.'s ring of life.
The exception being the trees depicted in Fig.~\ref{ring}
with the leaves marked by a filled-in square removed.
For  these 5 trees the associated phylogenetic
network ${\mathcal N}(T_1,\ldots, T_5)$ produced by the $MY$-closure
approach is depicted
in Fig.~\ref{lake-network}(left). 

In addition  and with the exception of one instance where one
split in $\Sigma(T_1,\ldots,T_5)$,
that is, the set of all splits displayed
by $T_1,\ldots,T_5$, was not extended to a full split
by our closure rules and thus was not displayed by
${\mathcal N}(T_1,\ldots, T_5)$ our rules always
generated a minimum collection of splits so that
${\mathcal N}(T_1,\ldots, T_5)$
displayed all the splits in $\Sigma(T_1,\ldots,T_5)$.

Interestingly, both the $Z$-closure super-network and Q-imputation approach
seemed to struggle with
this example with, in the case of $Z$-closure super-network,
either yielding a very complex network
${\mathcal N}(T_1,\ldots, T_5)$
in which numerous extensions of one and the same split in
$\Sigma(T_1,\ldots,T_5)$
was displayed (see Fig. \ref{lake-network}(right)) or
${\mathcal N}(T_1,\ldots, T_5)$
displayed only a  subset of splits in  $\Sigma(T_1,\ldots,T_5)$
(Q-imputation).\\

\noindent{\bf Acknowledgment}
The authors thank the referees for their helpful comments.

\end{document}